\documentclass[12pt,leqno]{article} 
\usepackage{amssymb} 
\usepackage{amsbsy}
\usepackage{amsmath} 
\usepackage{amsfonts} 
\usepackage{amscd}
\usepackage{graphicx,psfrag} 
\usepackage{color} 
 
\newcommand{\email}[1]{{\small E-mail: {\textsf {#1}}}}
\newcommand{\http}[1]{{\small Internet: {\textsf {#1}}}}
\newtheorem{thm}{Theorem} 
\newtheorem{prop}[thm]{Proposition} 
\newtheorem{lem}[thm]{Lemma} 
\newtheorem{cor}[thm]{Corollary} 
\newtheorem{rem}[thm]{Remark}

\newcommand{\beqn}{\begin{equation}} 
\newcommand{\eeqn}{\end{equation}} 
\newcommand{\bear}{\begin{eqnarray}} 
\newcommand{\eear}{\end{eqnarray}} 
\newcommand{\bean}{\begin{eqnarray*}} 
\newcommand{\eean}{\end{eqnarray*}} 
\newcommand\qed{\hfill $\square$} 
\newcommand{\RR}{\mathbb R}

\newcommand{\e}{\varepsilon} 
\begin{document} 

\title{Non-diffusive large time behaviour for\\ 
a degenerate viscous Hamilton-Jacobi equation}

\author{Philippe Lauren\c cot\footnote{Institut de 
Math\'ematiques de Toulouse, CNRS UMR~5219, Universit\'e de Toulouse, 118 route de Narbonne, F--31062 Toulouse Cedex 9, 
France. \email{Philippe.Laurencot@math.univ-toulouse.fr},
\http{http://www.mip.ups-tlse.fr/$\sim$laurenco/}}}
\date{\today}
\maketitle

\begin{abstract}
The convergence to non-diffusive self-similar solutions is investigated for non-negative solutions to the Cauchy problem $\partial_t u = \Delta_p u + \vert\nabla u\vert^q$ when the initial data converge to zero at infinity. Sufficient conditions on the exponents $p>2$ and $q>1$ are given that guarantee that the diffusion becomes negligible for large times and the $L^\infty$-norm of $u(t)$ converges to a positive value as $t\to\infty$.
\end{abstract}

\section{Introduction}\label{i}

\setcounter{thm}{0}
\setcounter{equation}{0}

The quasilinear degenerate parabolic equation
\beqn
\label{a1}
\partial_t u = \Delta_p u + \vert\nabla u\vert^q\ \,, 
\quad (t,x)\in Q_\infty:=(0,\infty)\times\RR^N\,,
\eeqn
includes two competing mechanisms acting on the space variable $x$, a degenerate diffusion $\Delta_p u$ involving the $p$-Laplacian operator defined by
$$
\Delta_p u := \mbox{ div } \left( \vert\nabla u\vert^{p-2}\ \nabla u
\right)\,, \quad p>2\,,
$$
and a source term $|\nabla u|^q$, $q>1$, depending solely on the gradient of $u$. The aim of this work is to identify a range of the parameters $p$ and $q$ for which the large time behaviour of non-negative solutions to \eqref{a1} is dominated by the source term. More precisely, we consider the Cauchy problem and supplement \eqref{a1} with the initial condition
\beqn
\label{a2}
u(0) = u_0\ge 0\,, \quad x\in\RR^N\,.
\eeqn
Throughout the paper, the initial condition $u_0$ is assumed to fulfill
\begin{equation}
\label{a3}
u_0 \in \mathcal{C}_0(\RR^N) \cap W^{1,\infty}(\RR^N)\,, \quad u_0\ge 0\,, \quad u_0\not\equiv 0\,, 
\end{equation}
where 
$$
\mathcal{C}_0(\RR^N) := \left\{ w \in \mathcal{BC}(\RR^N)\; :\; \lim_{R\to\infty} \sup_{\{|x|\ge R\}}{\{|w(x)|\}} = 0 \right\}\,,
$$
and $\mathcal{BC}(\RR^N):=\mathcal{C}(\RR^N)\cap L^\infty(\RR^N)$. 

For such an initial condition, the Cauchy problem \eqref{a1}, \eqref{a2} has a unique non-negative (viscosity) solution $u\in\mathcal{BC}([0,\infty)\times\RR^N)$ (see Proposition~\ref{prb1} below). Moreover, $t\longmapsto \|u(t)\|_\infty$ is a non-increasing function and has a limit $M_\infty\in [0,\|u_0\|_\infty]$ as $t\to\infty$. Our main result is then the following:

\begin{thm}\label{tha1}
Assume that $p>2$ and $q\in (1,p)$. Consider a non-negative function $u_0$ satisfying \eqref{a3} and let $u$ be the corresponding (viscosity) solution to \eqref{a1}, \eqref{a2}. Assume further that 
\beqn
\label{a4}
M_\infty := \lim_{t\to\infty} \Vert u(t)\Vert_\infty > 0\,.
\eeqn
Then
\beqn
\label{a5}
\lim_{t\to\infty} \Vert u(t) - h_\infty(t)\Vert_\infty = 0\,,
\eeqn
where $h_\infty$ is given by
\beqn
\label{a6}
h_\infty(t,x) := H_\infty\left( \frac{x}{t^{1/q}} \right) \quad\mbox{ and }\quad H_\infty(x) := \left( M_\infty - \gamma_q\ \vert x\vert^{q/(q-1)} \right)_+
\eeqn
for $(t,x)\in Q_\infty$ and $\gamma_q:= (q-1)\ q^{-q/(q-1)}$.
\end{thm} 
Here and below, $r_+:=\max{\{ r , 0 \}}$ denotes the positive part of the real number $r$.

The convergence \eqref{a5} clearly indicates that the large time behaviour of non-negative solutions to \eqref{a1}, \eqref{a2} fulfilling the condition \eqref{a4} is governed by the gradient source term. Indeed, $h_\infty$ is actually a self-similar solution to the Hamilton-Jacobi equation
\begin{equation}
\label{a7}
\partial_t h = \vert\nabla h\vert^q\ \,, \quad (t,x)\in Q_\infty\,,
\end{equation} 
and an alternative formula for $h_\infty$ reads
\beqn
\label{a8}
h_\infty(t,x) = \sup_{y\in\RR^N}\left\{ M_\infty\ \mathbf{1}_{\{0\}}(y)  - \gamma_q\ \frac{|x-y|^{q/(q-1)}}{t^{1/(q-1)}} \right\}
\eeqn
for $(t,x)\in [0,\infty)\times\RR^N$, $\mathbf{1}_{\{0\}}$ denoting the indicator function of the singleton set $\{0\}$. The formula \eqref{a8} is the well-known Hopf-Lax-Oleinik representation formula for viscosity solutions to \eqref{a7} (see, e.g., \cite[Chapter~3]{Ev98}) and $h_\infty$ turns out to be the unique viscosity solution in $\mathcal{BUC}(Q_\infty)$ to \eqref{a7} with the bounded and upper semicontinuous initial condition $h_\infty(0,x)=\mathbf{1}_{\{0\}}(x)$ for $x\in\RR^N$ \cite{St02}.

\begin{rem}\label{rea1}
The convergence \eqref{a5} also holds true for the viscosity solution to the Hamilton-Jacobi equation \eqref{a7} with a non-negative initial condition $u_0\in\mathcal{C}_0(\RR^N)$ but with $\|u_0\|_\infty$ instead of $M_\infty$ in the formula \eqref{a6} giving $H_\infty$. For \eqref{a1}, \eqref{a2}, the constant $M_\infty$ takes into account that, though negligible for large times, the diffusion erodes the supremum of $u$ during the time evolution.
\end{rem}

For $p=2$, Theorem~\ref{tha1} is also valid and is proved in \cite{BKL04}, the proof relying on a rescaling technique: The crucial step is then to identify the possible limits of the rescaled sequence and this is done by an extensive use of the Hopf-Lax-Oleinik representation formula. The proof we perform here is of a completely different nature and relies on the relaxed half-limits method introduced in \cite{BP88}. A similar approach has been used in \cite{NR99} and \cite{Ro01} to investigate the large time behaviour of solutions to first-order Hamilton-Jacobi equations $\partial_t w + H(x,\nabla w)=0$ in $Q_\infty$. It has also been used in \cite{LV07} to study the convergence to non-diffusive localized self-similar patterns for non-negative and compactly supported solutions to $\partial_t w -\Delta_p w + |\nabla w|^q=0$ in $Q_\infty$ when $p>2$ and $q\in (1,p-1)$. 

\medskip

In order to apply Theorem~\ref{tha1}, one should check whether there are non-negative solutions to \eqref{a1}, \eqref{a2} for which \eqref{a4} holds true. The next result provides sufficient conditions for \eqref{a4} to be fulfilled. 

\begin{thm}\label{tha2}
Assume that $p>2$ and $q>1$. Consider a non-negative function $u_0$ satisfying \eqref{a3} and let $u$ be the corresponding solution to \eqref{a1}, \eqref{a2}. Introducing
\beqn
\label{a9}
q_\star := p - \frac{N}{N+1}\,,
\eeqn
then $u$ fulfills \eqref{a4} if
\begin{itemize}
\item[(a)] either $q\in (1,q_\star]$,
\item[(b)] or $q\in (q_\star,p)$, $u_0\in W^{2,\infty}(\RR^N)$, and 
\beqn
\label{a10}
\Vert u_0\Vert_\infty > \kappa_0\ \left| \inf_{y\in\RR^N}{\left\{ \Delta_p u_0(y) \right\}} \right|^{(p-q)/q}\,.
\eeqn
for some $\kappa_0>0$ which depends only on $N$, $p$, and $q$. 
\end{itemize}
\end{thm}

A similar result is already available for $p=2$ and has been established in \cite{BKL04,Gi05}. The proof of Theorem~\ref{tha2} for $q\in (p-1,p)$ and $p>2$ borrows some steps from the case $p=2$. However, it relies on semiconvexity estimates for solutions to \eqref{a1}, \eqref{a2} which seem to be new for $p>2$ and $q\in (1,p)$ and are stated now.

\begin{prop}\label{pra3}
Assume that $p>2$ and $q\in (1,p]$. Let $u$ be the viscosity solution to \eqref{a1}, \eqref{a2} with initial condition $u_0\in\mathcal{BUC}(\RR^N)$ (that is, $u_0\in\mathcal{BC}(\RR^N)$ and is uniformly continuous in $\RR^N$). Then $\nabla u(t)$ belongs to $L^\infty(\RR^N)$ for each $t>0$ and there is $\kappa_1>0$ depending only on $N$, $p$, and $q$ such that
\beqn
\label{a11}
\Delta_p u(t,x) \ge - \kappa_1\ \|u(s)\|_\infty^{(p-q)/q}\ (t-s)^{-p/q}\,, \quad t>s\ge 0\,,
\eeqn
in the sense of distributions. In addition, if $u_0\in W^{1,\infty}(\RR^N)$, there holds
\beqn
\label{a12}
\Delta_p u(t,x) \ge - \frac{N (p-1)}{q (q-1)}\ \frac{\|\nabla u_0\|_\infty^{p-q}}{t}
\eeqn
for $t>0$ in the sense of distributions.
\end{prop}

The proof of Proposition~\ref{pra3} relies on the comparison principle combined with a gradient estimate established in \cite{BtL08}. 

Similar semiconvexity estimates for solutions to \eqref{a1}, \eqref{a2} have already been obtained in \cite{Ha93} and \cite[Lemma~5.1]{LT01} for $p=q=2$, in \cite[Proposition~3.2]{BKL04} for $p=2$ and $q\in (1,2]$, and in \cite[Theorem~1]{EM94} for $p=q>2$. We extend these results to the range $p>2$ and $q\in (1,p]$. As we shall see below, the estimate \eqref{a11} plays an important role in the proof of Theorem~\ref{tha2} and is also helpful to construct a subsolution in the proof of Theorem~\ref{tha1}. 

Let us finally emphasize that the validity of Proposition~\ref{pra3} is not restricted to non-negative solutions and that the solutions to the Hamilton-Jacobi equation \eqref{a7} also enjoy the semiconvexity estimates \eqref{a11} and \eqref{a12}. These two estimates thus stem from the reaction term $|\nabla u|^q$ and not from the diffusion.

\medskip

In the next section, we recall the well-posedness of \eqref{a1}, \eqref{a2} in $\mathcal{BUC}(\RR^N)$, as well as some properties of the solutions established in \cite{BtL08}. We also show the finite speed of propagation of the support for non-negative compactly supported initial data. Section~\ref{s} is devoted to the proof of the semiconvexity estimates (Proposition~\ref{pra3}) and Section~\ref{css} to that of Theorem~\ref{tha1}. Theorem~\ref{tha2} is shown in the last section, its proof combining arguments of \cite{BKL04,Gi05,LS05} used to established analogous results when $p=2$. 

\medskip

Throughout the paper, $C$ and $C_i$, $i\ge 1$, denote positive constants depending only on $p$, $q$, and $N$. Dependence upon additional parameters will be indicated explicitly. Also, $\mathcal{M}_N(\RR)$ denotes the space of real-valued $N\times N$ matrices and $\delta_{ij}=1$ if $i=j$ and $\delta_{ij}=0$ if $i\ne j$, $1\le i,j\le N$. Given a matrix $A=(a_{ij})\in \mathcal{M}_N(\RR)$, $tr(A)$ denotes its trace and is given by $tr(A):=\sum a_{ii}$. 

\section{Preliminary results}\label{p}

\setcounter{thm}{0}
\setcounter{equation}{0}

Let us first recall the well-posedness (in the framework of viscosity solutions) of \eqref{a1}, \eqref{a2}, together with some properties of the solutions established in \cite{BtL08}.

\begin{prop}\label{prb1}
Consider a non-negative initial condition $u_0\in\mathcal{BUC}(\RR^N)$. There is a unique non-negative viscosity solution
$u \in \mathcal{BC}([0,\infty)\times \RR^N)$ to \eqref{a1}, \eqref{a2} such that 
\beqn
\label{b1}
0\le u(t,x)\le \|u_0\|_\infty\,, \quad (t,x)\in Q_\infty\,,
\eeqn 
\beqn
\label{b2}
\left\| \nabla u(t) \right\|_\infty \le \min{\left\{ C_1\ \Vert u(s)\Vert_\infty^{1/q}\ (t-s)^{-1/q} , \|\nabla u(s)\|_\infty \right\}}\,,
\eeqn
and
\beqn
\label{b3}
\int_{\RR^N} (u(t,x)-u(s,x))\ \vartheta(x)\ dx + \int_s^t \int_{\RR^N}
\left( |\nabla u|^{p-2} \nabla u \cdot \nabla\vartheta - |\nabla u|^q\
\vartheta \right)\ dxd\tau = 0
\eeqn
for $t>s\ge 0$ and $\vartheta\in\mathcal{C}_0^\infty(\RR^N)$. 
In addition, $t\longmapsto \|u(t)\|_\infty$ is a non-increasing function.
\end{prop}

\noindent\textbf{Proof.} We put $\tilde{u}_0:= \|u_0\|_\infty - u_0$. As $\tilde{u}_0$ is a non-negative function in $\mathcal{BUC}(\RR^N)$, it follows from \cite[Theorem~1.1]{BtL08} that there is a unique non-negative viscosity solution $\tilde{u}$ to 
\beqn
\label{b4}
\partial_t \tilde{u} - \Delta_p \tilde{u} + \vert\nabla \tilde{u}\vert^q = 0\ \,, 
\quad (t,x)\in Q_\infty:=(0,\infty)\times\RR^N\,,
\eeqn
with initial condition $\tilde{u}(0,x)=\tilde{u}_0(x)$ for $x\in\RR^N$. It also satisfies $0\le \tilde{u}(t,x)\le \|u_0\|_\infty$ and 
$$
\int_{\RR^N} (\tilde{u}(t,x)-\tilde{u}(s,x))\ \vartheta(x)\ dx + \int_s^t \int_{\RR^N} \left( |\nabla \tilde{u}|^{p-2} \nabla \tilde{u} \cdot \nabla\vartheta + |\nabla \tilde{u}|^q\ \vartheta \right)\ dxd\tau = 0
$$
for $t>s\ge 0$, $x\in\RR^N$, and $\vartheta\in\mathcal{C}_0^\infty(\RR^N)$. In addition, $\nabla\tilde{u}(t)$ belongs to $L^\infty(\RR^N)$ for each $t>0$ and
$$
\|\nabla\tilde{u}(t)\|_\infty \le C_1\ \left\|\tilde{u}_0\right\|_\infty^{1/q}\ t^{-1/q}
$$
by \cite[Lemma~4.1]{BtL08}. Setting $u:=\|u_0\|_\infty-\tilde{u}$, we readily deduce from the properties of $\tilde{u}$ that $u$ is a non-negative viscosity solution to \eqref{a1}, \eqref{a2} satisfying \eqref{b1} and \eqref{b3}. Also, $\nabla u(t)$ belongs to $L^\infty(\RR^N)$ for each $t>0$. The uniqueness and the time monotonicity of $\|u\|_\infty$ then both follow from the comparison principle, see \cite{CIL92} or  \cite[Theorem~2.1]{GGIS91}. Finally, given $s\ge 0$, $(t,x)\mapsto \|u(s)\|_\infty -u(t+s,x)$ is the unique non-negative viscosity solution to the Cauchy problem \eqref{b4} with initial condition $x\mapsto \|u(s)\|_\infty-u(s,x)$ and we infer from \cite[Lemma~4.1]{BtL08} that
$$
\|\nabla u(t+s)\|_\infty \le C_1\ \left\| \|u(s)\|_\infty - u(s) \right\|_\infty^{1/q}\ t^{-1/q} \le C_1\ \| u(s)\|_\infty^{1/q}\ t^{-1/q}
$$
for $t>0$, whence \eqref{b2}. \qed

We next turn to the propagation of the support of non-negative solutions to \eqref{a1}, \eqref{a2} with non-negative compactly supported initial data. 

\begin{prop}\label{prb2}
Consider a non-negative solution $u$ to \eqref{a1}, \eqref{a2} with an initial condition $u_0$ satisfying \eqref{a3}. Assume further that $u_0$ is compactly supported in a ball $B(0,R_0)$ of $\RR^N$ for some $R_0>0$. Then $u(t)$ is compactly supported for each $t\ge 0$.
\end{prop}

\noindent\textbf{Proof.} We argue by comparison with travelling wave solutions. By \cite[Application~9.4]{GK04}, there is a travelling wave solution $w$ to the convection-diffusion equation
\begin{equation}
\label{b5}
\partial_t w - \partial_1^2\left( w^{p-1} \right) + \partial_1\left( w^q \right) = 0\,, \quad (t,x_1)\in (0,\infty)\times\RR\,,
\end{equation}
with wave speed unity. It is given by $w(t,x_1)=f(x_1-t)$ for $(t,x_1)\in (0,\infty)\times\RR$, the function $f$ being implicitly defined by 
$$
(p-1)\ \int_0^{f(y)} \frac{z^{p-3}}{1-z^{q-1}}\ dz = (-y)_+\,, \quad y\in\RR\,.
$$
In particular, $f$ satisfies $f(y)=0$ if $y>0$ and $f(y)\to 1$ as $y\to -\infty$. Introducing
$$
F(y) := \int_y^\infty f(z)\ dz\,, \quad y\in\RR\,,
$$
the properties of $f$ ensure that $F$ is a decreasing function on $(-\infty,0)$ with $F(y)=0$ if $y>0$, $F(y)\le |y|$ if $y<0$, and $F(y)\to\infty$ as $y\to -\infty$. There is therefore a unique $\mu\in (-\infty,0)$ such that $F(R_0+\mu)=\|u_0\|_\infty$. In addition, it readily follows from \eqref{b5} and the invariance by translation of \eqref{a1} that $W_\mu(t,x):=F(x_1+\mu-t)$ is a travelling wave solution to \eqref{a1}. Now, $u$ and $W_\mu$ are both solutions to \eqref{a1} in $(0,\infty)\times H_+$, the half-space $H_+$ being defined by $H_+:=\left\{ x\in\RR^N\; : \; x_1>R_0\right\}$. Owing to the monotonicity of $F$, the bound $0\le f\le 1$, and \eqref{b1}, we have also
$$
u_0(x) - W_\mu(0,y) = 0 - W_\mu(0,y) \le W_\mu(0,x) - W_\mu(0,y) \le |x-y|
$$
for $x\in H_+$ and $y\in H_+$,
\bean
u(t,x) - W_\mu(t,y) & \le & \|u_0\|_\infty - W_\mu(t,y) \\
& \le & F(R_0+\mu-t) - W_\mu(t,y) = W_\mu(t,x) - W_\mu(t,y) \le |x-y|
\eean
for $t>0$, $x\in\partial H_+$, $y\in H_+$, and 
$$
u(t,x) - W_\mu(t,y) \le \|u_0\|_\infty - F(R_0+\mu-t) \le 0
$$
for $t>0$, $x\in H_+$, $y\in\partial H_+$. We are then in a position to use the comparison principle stated in \cite[Theorem~2.1]{GGIS91} to conclude that $u(t,x)\le W_\mu(t,x)$ for $(t,x)\in (0,\infty)\times H_+$. Consequently,  $u(t,x)\le F(x_1+\mu-t)=0$ if $t\ge 0$ and $x_1>\max{\{R_0,t-\mu\}}$, and the rotational invariance of \eqref{a1} allows us to conclude that $u(t,x)=0$ for $t\ge 0$ and $|x|>\max{\{R_0,t-\mu\}}$. \qed

\medskip

We finally recall the convergence to self-similar solutions for non-negative and compactly supported solutions to the $p$-Laplacian equation \cite{KV88}
\beqn
\label{b6}
\partial_t \varphi = \Delta_p \varphi\ \,, \quad (t,x)\in Q_\infty\,.
\eeqn

\begin{prop}\label{prb3}
Let $\varphi_0$ be a non-negative and compactly supported function in $L^1(\RR^N)$ and $\varphi$ denote the unique weak solution to \eqref{b6} with initial condition $\varphi_0$. Then
\beqn
\label{b7}
\lim_{t\to\infty} t^{(N(r-1))/(r(N(p-2)+p))}\ \left\| \varphi(t) - \mathcal{B}_{\|\varphi_0\|_1}(t) \right\|_r = 0 \;\;\mbox{ for }\;\; r\in [1,\infty]\,,
\eeqn
where $\mathcal{B}_L$ denotes the Barenblatt solution to \eqref{b6} given by
\bean
\mathcal{B}_L(t,x) & := & t^{-N/(N(p-2)+p)}\ b_L\left( x t^{-1/(N(p-2)+p)}\right)\,, \\
b_L(x) & := & \left( C_2\ L^{(p(p-2))/((p-1)(N(p-2)+p))} - C_3\ |x|^{p/(p-1)} \right)_+^{(p-1)/(p-2)}
\eean
for $(t,x)\in (0,\infty)\times\RR^N$ and $L>0$.
\end{prop}

The convergence \eqref{b7} is proved in \cite[Theorem~2]{KV88} for $r=\infty$. As $\varphi_0$ is compactly supported, so is $\varphi(t)$ for each $t>0$ and the support of $\varphi(t)$ is included in $B\left( 0,C_4(\varphi_0)\ t^{1/(N(p-2)+p)} \right)$ for $t\ge 1$ \cite[Proposition~2.2]{KV88}. Combining this property with \cite[Theorem~2]{KV88}  readily provide the convergence \eqref{b7} for all $r\in [1,\infty)$. 

\section{Semiconvexity}\label{s}

\setcounter{thm}{0}
\setcounter{equation}{0}

In this section, we prove the semiconvexity estimates \eqref{a11} and \eqref{a12}. To this end, we would like to derive an equation for $\Delta_p u$ to which we could apply the comparison principle. The poor regularity of $u$ however does not allow to perform directly such a computation and an approximation procedure is needed. As a first step, we report the following result:

\begin{lem}\label{lec1}
Let $a$ and $b$ be two non-negative function in $\mathcal{C}^\infty([0,\infty))$ satisfying 
\bear
\label{c0a}
& & a(r)>0\,, \quad a'(r)>0\,, \quad a'(r)\ b'(r) - a(r)\ b''(r)>0\,,\\
\label{c0b}
& & c(r) := 2\ \left( \frac{b'}{a} \right)(r) + \frac{4r\ (a\ b''-a'\ b')(r)}{a^2(r) + 2 r\ a(r)\ a'(r)}\ge 0\,.
\eear
Consider a classical solution $v$ to 
\beqn
\label{c1}
\partial_t v - \mbox{ div }\left( a\left( |\nabla v|^2 \right)\ \nabla v \right) = b\left( |\nabla v|^2 \right)\,, \quad (t,x)\in Q_\infty\,,
\eeqn
and put 
$$
w:= \mbox{ div }\left( a\left( |\nabla v|^2\right)\ \nabla v \right) \;\;\mbox{ and }\;\; z_i := a\left( |\nabla v|^2\right)\ \partial_i v
$$
for $i\in\{1,\ldots,N\}$. Then
\beqn
\label{c2}
\partial_t w - \mathcal{L} w - \mathcal{V}\cdot \nabla w - \frac{c\left( |\nabla v|^2 \right)}{N}\ w^2 \ge 0 \;\;\mbox{ in }\;\; Q_\infty\,,
\eeqn
where
\bean
\mathcal{L} w & := & \sum_{i,j} \partial_i \left( a\left( |\nabla v|^2 \right)\ E_{ij}\ \partial_j w \right)\,, \quad \mathcal{V} := 2\ b'\left( |\nabla v|^2 \right)\ \nabla v\,, \\
E_{ij} & := & \delta_{ij} + 2\ \frac{a'}{a}\left( |\nabla v|^2\right)\ \partial_i v\ \partial_j v\,, \quad 1 \le i,j\le N\,.
\eean
\end{lem}

The proof of Lemma~\ref{lec1} borrows some steps from the proof of \cite[Theorem~1]{EM94} for $p=q>2$ but requires additional arguments to handle the term coming from the fact that $q\ne p$. In particular, we recall the following elementary result which will be helpful to estimate this term.

\begin{lem}\label{lec2}
Let $A$ and $B$ be two symmetric matrices in $\mathcal{M}_N(\RR)$ and put $M:=ABA$. Then $M$ is a symmetric matrix in $\mathcal{M}_N(\RR)$ and
\beqn
\label{c3}
|MX|^2 \le tr\left( M^2 \right)\ |X|^2 \;\;\mbox{ for }\;\; X\in\RR^N\,.
\eeqn
\end{lem}

\noindent\textbf{Proof of Lemma~\ref{lec1}.} We first note that
\bear
\label{c4}
\partial_j z_i & = & a\left( |\nabla v|^2\right)\ \sum_k E_{ik}\ \partial_k \partial_j v\,,\\
\label{c5}
\partial_t z_i & = & a\left( |\nabla v|^2\right)\ \sum_k E_{ik}\ \partial_k \partial_t v\,,
\eear
for $1\le i,j\le N$. According to the definition of $w$, we infer from \eqref{c1}, \eqref{c4}, and \eqref{c5} that
\bean
\partial_t w & = & \sum_{i,k} \partial_i \left( a\left( |\nabla v|^2\right)\ E_{ik}\ \partial_k \partial_t v \right) \\
& = &  \sum_{i,k} \partial_i \left( a\left( |\nabla v|^2\right)\ E_{ik}\ \partial_k \left( w + b\left( |\nabla v|^2\right) \right) \right) \\
& = & \mathcal{L}w + 2\ \sum_{i,k} \partial_i\left( (ab')\left( |\nabla v|^2\right)\ E_{ik}\ \sum_j \partial_j v\ \partial_j \partial_k v \right) \\
& = & \mathcal{L}w + 2\ \sum_{i,j} \partial_i\left( (ab')\left( |\nabla v|^2\right)\ \partial_j v\ \sum_k E_{ik}\  \partial_k \partial_j v \right) \\
& = & \mathcal{L}w + 2\ \sum_{i,j} \partial_i\left( \left( \frac{b'}{a} \right)\left( |\nabla v|^2\right)\ z_j\ \partial_j z_i \right) \\
& = & \mathcal{L}w + 4\ \sum_{i,j} \left( \frac{b'}{a} \right)'\left( |\nabla v|^2\right)\ \sum_k \partial_k v\ \partial_k \partial_i v\ z_j\ \partial_j z_i \\ 
& + & 2\ \sum_{i,j} \left( \frac{b'}{a} \right)\left( |\nabla v|^2\right)\ \partial_i z_j\ \partial_j z_i + 2\ \sum_{i,j} \left( \frac{b'}{a} \right)\left( |\nabla v|^2\right)\ z_j\ \partial_j \partial_i z_i\,.
\eean
Since $w=\sum \partial_i z_i$, the last term of the right-hand side of the above inequality is equal to $\mathcal{V}\cdot\nabla w$ and 
\bear
\label{c6}
\partial_t w & = & \mathcal{L}w + \mathcal{V}\cdot\nabla w + 4\ \left[ a\ \left( \frac{b'}{a} \right)' \right]\left( |\nabla v|^2\right)\ \sum_{i,j,k} \partial_j v\ \partial_k v\ \partial_k \partial_i v\ \partial_j z_i \\ 
\nonumber
& + & 2\ \left( \frac{b'}{a} \right)\left( |\nabla v|^2\right)\ \sum_{i,j}  \partial_i z_j\ \partial_j z_i \,.
\eear
On the one hand, introducing the matrix $\mathcal{E}:=(E_{ij})$ and the Hessian matrix $D^2 v=(\partial_i \partial_j v)$ of $v$, we infer from \eqref{c4} that 
\bear
\nonumber
\sum_{i,j}  \partial_i z_j\ \partial_j z_i & = & a^2\left( |\nabla v|^2 \right)\ \sum_{i,j,k,l} E_{ik}\ \partial_k \partial_j v\ E_{jl}\ \partial_l \partial_i v \\
\nonumber
& = & a^2\left( |\nabla v|^2 \right)\ \sum_{i,j} \left( \mathcal{E}\ D^2 v \right)_{ij}\ \left( \mathcal{E}\ D^2 v \right)_{ji}\\
\label{c7}
\sum_{i,j}  \partial_i z_j\ \partial_j z_i & = & a^2\left( |\nabla v|^2 \right)\ tr\left( \left( \mathcal{E}\ D^2 v \right)^2\right)\,.
\eear
On the other hand, using once more \eqref{c4}, we obtain
\bear
\nonumber
\sum_{i,j,k} \partial_j v\ \partial_k v\ \partial_k \partial_i v\ \partial_j z_i & = & a\left( |\nabla v|^2 \right)\ \sum_{i,j,k,l} \partial_j v\ \partial_k v\ \partial_k \partial_i v\ E_{il}\ \partial_l \partial_j v \\
\nonumber
& = & a\left( |\nabla v|^2 \right)\ \sum_{i,l} \left( \sum_k \partial_i \partial_k v\ \partial_k v \right)\ E_{il}\ \left( \sum_j \partial_l \partial_j v\ \partial_j v \right) \\
\label{c8}
\quad \sum_{i,j,k} \partial_j v\ \partial_k v\ \partial_k \partial_i v\ \partial_j z_i & = & a\left( |\nabla v|^2 \right)\ \left\langle D^2 v\ \nabla v , (\mathcal{E}\ D^2 v)\ \nabla v \right\rangle\,.
\eear
Inserting \eqref{c7} and \eqref{c8} in \eqref{c6}, we end up with
\bear
\label{c9}
\partial_t w & = & \mathcal{L}w + \mathcal{V}\cdot\nabla w + 2\ \left( ab' \right)\left( |\nabla v|^2\right)\ tr\left( \left( \mathcal{E}\ D^2 v \right)^2\right) \\
\nonumber
& + & 4\ \left( a\ b'' - a'\ b' \right)\left( |\nabla v|^2\right)\ \left\langle D^2 v\ \nabla v , (\mathcal{E}\ D^2 v)\ \nabla v \right\rangle\,.
\eear

We next observe that
\beqn
\label{c10}
\mathcal{E}\ \nabla v = \left( 1 + 2\ |\nabla v|^2\ \left( \frac{a'}{a} \right)\left( |\nabla v|^2\right)\right)\ \nabla v
\eeqn
and that, for $X\in\RR^N$, 
$$
\left\langle \mathcal{E}\ X , X \right\rangle = |X|^2 + 2\ \left( \frac{a'}{a} \right)\left( |\nabla v|^2\right)\ \langle X , \nabla v \rangle^2 \ge |X|^2
$$
as $a$ and $a'$ are both positive by \eqref{c0a}. Consequently, $\mathcal{E}$ is a positive definite symmetric matrix in $\mathcal{M}_N(\RR)$ and there exists a positive definite matrix $\mathcal{E}_{1/2}$ such that $\mathcal{E}_{1/2}^2=\mathcal{E}$. We then infer from the definition of $\mathcal{E}_{1/2}$, \eqref{c10}, and Lemma~\ref{lec2} (with $A=\mathcal{E}_{1/2}$, $B=D^2 v$ and $X=\mathcal{E}_{1/2}^{-1}\ \nabla v$) that 
\bean
\left\langle D^2 v\ \nabla v , (\mathcal{E}\ D^2 v)\ \nabla v \right\rangle & = & \left| (\mathcal{E}_{1/2}\ D^2 v)\ \nabla v \right|^2 \\
& = & \left| \left( \left( \mathcal{E}_{1/2}\ D^2 v\ \mathcal{E}_{1/2} \right)\ \mathcal{E}_{1/2}^{-1} \right)\ \nabla v \right|^2 \\
& \le & tr\left( \mathcal{E}_{1/2}\ D^2 v\ \mathcal{E}_{1/2}\ \mathcal{E}_{1/2}\ D^2 v\ \mathcal{E}_{1/2} \right)\ \left\langle \mathcal{E}_{1/2}^{-1}\ \nabla v ,  \mathcal{E}_{1/2}^{-1}\ \nabla v \right\rangle \\
& \le & tr\left( \left( \mathcal{E}\ D^2 v \right)^2\right)\ \left\langle \nabla v ,  \mathcal{E}^{-1}\ \nabla v \right\rangle \\
& \le & tr\left( \left( \mathcal{E}\ D^2 v \right)^2\right)\ |\nabla v|^2\ \left( 1 + 2\ |\nabla v|^2\ \left( \frac{a'}{a} \right)\left( |\nabla v|^2\right)\right)^{-1}\,.
\eean
Owing to the non-positivity \eqref{c0a} of $a\ b'' - a'\ b'$, we deduce from \eqref{c9} and the above inequality that
$$
\partial_t w \ge \mathcal{L}w + \mathcal{V}\cdot\nabla w + \left( a^2\ c\right)\left( |\nabla v|^2\right)\ tr\left( \left( \mathcal{E}\ D^2 v \right)^2\right)\,,
$$
the function $c$ being defined in \eqref{c0b}. We finally use the inequality
$$
tr\left( A^2 \right) \ge \frac{1}{N}\ tr(A)^2\,, \quad A\in\mathcal{M}_N(\RR)\,,
$$
the identity
$$
w = \sum_i \partial_i z_i = a\left( |\nabla v|^2 \right)\ tr\left( \mathcal{E}\ D^2 v \right)\,,
$$
and the non-negativity \eqref{c0b} of $c$ to conclude that
\bean
\partial_t w & \ge & \mathcal{L}w + \mathcal{V}\cdot\nabla w + \frac{1}{N}\ \left( a^2\ c\right)\left( |\nabla v|^2\right)\ tr\left( \mathcal{E}\ D^2 v \right)^2 \\
& \ge & \mathcal{L}w + \mathcal{V}\cdot\nabla w + \frac{c\left( |\nabla v|^2\right)}{N}\ w^2\,,
\eean
and complete the proof. \qed

\medskip

\noindent\textbf{Proof of Proposition~\ref{pra3}.} To be able to use Lemma~\ref{lec1}, we shall first construct a suitable approximation of \eqref{a1}, \eqref{a2}. Such a construction has already been performed in \cite{BtL08} for similar purposes and we recall it now. Given $u_0$ satisfying \eqref{a3}, there is a sequence of functions $(u_{0,k})_{k\ge 1}$ such that, for each integer $k\ge 1$, $u_{0,k}\in \mathcal{BC}^\infty(\RR^N)$, $u_0 \le u_{0,k+1} \le u_{0,k}$, and $(u_{0,k},\nabla u_{0,k})_k$ converge towards $(u_0,\nabla u_0)$ uniformly on every compact subset of $\RR^N$ as $k\to\infty$. Next, for $\varepsilon\in (0,1)$ and $r\ge 0$, we set
$$
a_\varepsilon(r) := \left( r + \varepsilon^2 \right)^{(p-2)/2} \;\;\mbox{ and }\;\; b_\varepsilon(r):=\left( r + \varepsilon^2 \right)^{q/2} - \varepsilon^q\,.
$$ 
Then the Cauchy problem
\bear
\label{c11}
\partial_t u_{k,\varepsilon} & = & \mbox{ div }\left( a_\varepsilon\left( |\nabla u_{k,\varepsilon}|^2 \right)\ \nabla u_{k,\varepsilon} \right) + b_\varepsilon\left( |\nabla u_{k,\varepsilon}|^2 \right) \,, \quad (t,x)\in Q_\infty\,, \\
\label{c12}
u_{k,\varepsilon}(0) & = & u_{0,k} + \varepsilon^\nu\,, \quad x\in\RR^N\,,
\eear
has a unique classical solution $u_{k,\varepsilon}$, the parameter $\nu>0$ depending $p$, $q$, and $N$ and being appropriately chosen. Furthermore,
\bear
\label{c13}
\|\nabla u_{k,\varepsilon}(t) \|_\infty & \le & \|\nabla u_{0,k}\|_\infty\,, \quad t\ge 0\,,\\
\label{c14}
\lim_{k\to\infty}\ \lim_{\varepsilon\to 0} u_{k,\varepsilon}(t,x) & = & u(t,x) \,,
\eear
the latter convergence being uniform on every compact subset of $[0,\infty)\times\RR^N$, see \cite[Section~3]{BtL08} (after performing the same change of unknown function as in the proof of Proposition~\ref{prb1}). 

Introducing
$$
c_\varepsilon(r) = 2\ \left( \frac{b_\varepsilon'}{a_\varepsilon} \right)(r) + \frac{4r\ (a_\varepsilon\ b_\varepsilon''-a_\varepsilon'\ b_\varepsilon')(r)}{a_\varepsilon^2(r) + 2 r\ a_\varepsilon(r)\ a_\varepsilon'(r)}\,, \quad r\ge 0\,,
$$
let us check that $a_\varepsilon$ and $b_\varepsilon$ fulfill the conditions \eqref{c0a} and \eqref{c0b}. Clearly, $a_\varepsilon>0$ and $a_\varepsilon'>0$ as $p>2$. Next, since $1<q\le p$, 
\bean
\left( a_\varepsilon'\ b_\varepsilon' - a_\varepsilon\ b_\varepsilon'' \right)(r) & = & \frac{q\ (p-q)}{4}\ \left( r+\varepsilon^2 \right)^{(p+q-6)/2} \ge 0\,, \\
c_\varepsilon(r) &= & q\ \frac{r(q-1)+\varepsilon^2}{r(p-1)+\varepsilon^2}\ \left( r+\varepsilon^2 \right)^{(q-p)/2} \ge 0\,.
\eean
We may then apply Lemma~\ref{lec1} to deduce that $w_{k,\varepsilon}:=\mbox{ div }\left( a_\varepsilon\left( |\nabla u_{k,\varepsilon}|^2 \right)\ \nabla u_{k,\varepsilon} \right)$ satisfies
$$
 \partial_t w_{k,\varepsilon} - \mathcal{L}_{k,\varepsilon} w_{k,\varepsilon} - \mathcal{V}_{k,\varepsilon}\cdot \nabla w_{k,\varepsilon} - \frac{c_\varepsilon\left( |\nabla u_{k,\varepsilon}|^2 \right)}{N}\ w_{k,\varepsilon}^2 \ge 0
$$
in $Q_\infty$. Observe next that the condition $1<q\le p$ implies that $c_\varepsilon$ is a non-increasing function. It then follows from \eqref{c13} that $c_\varepsilon\left( |\nabla u_{k,\varepsilon}|^2 \right) \ge c_\varepsilon\left( \|\nabla u_{0,k}\|_\infty^2 \right)$ and we end up with
\beqn
\label{c15}
 \partial_t w_{k,\varepsilon} - \mathcal{L}_{k,\varepsilon} w_{k,\varepsilon} - \mathcal{V}_{k,\varepsilon}\cdot \nabla w_{k,\varepsilon} - \frac{c_\varepsilon\left( \|\nabla u_{0,k}\|_\infty^2 \right)}{N}\ w_{k,\varepsilon}^2 \ge 0
\eeqn
in $Q_\infty$. Clearly, $t\longmapsto - N/\left( c_\varepsilon\left( \|\nabla u_{0,k}\|_\infty^2 \right)\ t \right)$ is a subsolution to \eqref{c15} and the comparison principle warrants that 
\beqn
\label{c15b}
w_{k,\varepsilon}(t,x) \ge - \frac{N}{c_\varepsilon\left( \|\nabla u_{0,k}\|_\infty^2 \right)\ t}\,, \quad (t,x)\in Q_\infty\,.
\eeqn
Letting $\varepsilon\to 0$ and $k\to\infty$ in the previous inequality with the help of \eqref{c14} gives \eqref{a12}. 

Next, since \eqref{a1} is autonomous, we infer from \eqref{b2} (with $s=0$) and \eqref{a12} that  
\bean
\Delta_p u(t,x) & \ge & - \frac{2N (p-1)}{q (q-1)}\ \frac{\|\nabla u(t/2)\|_\infty^{p-q}}{t} \\
& \ge & - \frac{2^{p/q} N (p-1)}{q (q-1)}\ C_1^{p-q}\ \| u_0\|_\infty^{(p-q)/q}\ t^{-p/q}\,,
\eean
whence \eqref{a11} for $s=0$. To prove the general case $s\in (0,t)$, we use again the fact that \eqref{a1} is autonomous. \qed

\medskip

We have a similar result when $u_0$ is more regular.

\begin{cor}\label{coc3}
Assume that $p>2$ and $q\in (1,p]$. Let $u$ be the solution to \eqref{a1}, \eqref{a2} with an initial condition $u_0$ satisfying $u_0\in W^{2,\infty}(\RR^N)$ in addition to \eqref{a3}. Then
\beqn
\label{c16}
\Delta_p u(t,x) \ge - \left| \inf_{y\in\RR^N} {\Delta_p u_0(y)} \right|
\eeqn
in the sense of distributions.
\end{cor}

\noindent\textbf{Proof.} Keeping the notations introduced in the proof of Proposition~\ref{pra3}, we readily infer from \eqref{c15} and the comparison principle that 
\beqn
\label{c17}
w_{k,\varepsilon}(t,x) \ge - \left| \inf_{y\in\RR^N} {\Delta_p u_{0,k}(y)} \right|\,, \quad (t,x)\in Q_\infty\,.
\eeqn
Owing to the regularity of $u_0$, it is possible to construct the sequence $(u_{0,k})_k$ such that it satisfies
$$
\lim_{k\to\infty} \inf_{y\in\RR^N} {\Delta_p u_{0,k}(y)} = \inf_{y\in\RR^N} {\Delta_p u_0(y)}\,.
$$
We may then pass to the limit first as $\varepsilon\to 0$ and then as $k\to\infty$ in \eqref{c17} and use \eqref{c14} and the above convergence to complete the proof. \qed

\medskip

Another useful consequence of the semiconvexity estimates derived in Proposition~\ref{pra3} is that the solution $u$ to \eqref{a1}, \eqref{a2} is a supersolution to a first-order Hamilton-Jacobi equation.

\begin{cor}\label{coc4}
Consider an initial condition $u_0$ satisfying \eqref{a3}. Setting $F(t,\xi_0,\xi):= \xi_0 - |\xi|^q + \kappa_1\ \|u_0\|_\infty^{(p-q)/q}\ t^{-p/q}$ for $t\in (0,\infty)$, $\xi_0\in\RR$, and $\xi\in\RR^N$ (recall that $\kappa_1$ is defined in \eqref{a11}), the solution $u$ to \eqref{a1}, \eqref{a2} is a supersolution to $F(t,\partial_t w,\nabla w)=0$ in $Q_\infty$.
\end{cor}

\noindent\textbf{Proof.} We still use the notations introduced in the proof of Proposition~\ref{pra3}. As $w_{k,\varepsilon}=\mbox{ div }\left( a_\varepsilon\left( |\nabla u_{k,\varepsilon}|^2 \right)\ \nabla u_{k,\varepsilon} \right)$, we infer from \eqref{c11} and \eqref{c15b} that
$$
\partial_t u_{k,\varepsilon} - b_\varepsilon\left( |\nabla u_{k,\varepsilon}|^2 \right) \ge - \frac{N}{c_\varepsilon\left( \|\nabla u_{0,k}\|_\infty^2 \right)\ t}
$$
in $Q_\infty$. We then use \eqref{c14} and the stability of viscosity solutions \cite{BdCD97,Bl94,CIL92} to pass to the limit as $\varepsilon\to 0$ and $k\to\infty$ in the previous inequality and conclude that $u$ is a supersolution to
$$
\partial_t w - |\nabla w|^q + \frac{N (p-1)}{q (q-1)}\ \frac{\|\nabla u_0\|_\infty^{p-q}}{t} = 0 \;\;\mbox{ in }\;\; Q_\infty\,.
$$

Now, fix $T\ge 0$. As \eqref{a1} is an autonomous equation, the function $(t,x)\longmapsto u(t+T,x)$ is the solution to \eqref{a1} with initial condition $u(T)$ and the above analysis allows us to conclude that $u$ is a supersolution to
$$
\partial_t w - |\nabla w|^q + \frac{N (p-1)}{q (q-1)}\ \frac{\|\nabla u(T)\|_\infty^{p-q}}{t-T} = 0 \;\;\mbox{ in }\;\; (T,\infty)\times\RR^N\,.
$$
We then use \eqref{b2} (with $T=t/2$) to complete the proof. \qed

\section{Convergence to self-similarity}\label{css}

\setcounter{thm}{0}
\setcounter{equation}{0}

We change the variables and the unknown function so that the convergence (\ref{a5}) is transformed to the convergence towards a steady state. More precisely, we introduce the self-similar (or scaling) variables 
$$
\tau = \frac{1}{q}\ \log{(1+t)}\,, \qquad y = \frac{x}{(1+t)^{1/q}}\,,
$$
and the new unknown function $v$ defined by
\beqn
\label{d1}
u(t,x) = v\left( \frac{\log{(1+t)}}{q} , \frac{x}{(1+t)^{1/q}} \right)\,, \qquad (t,x)\in [0,\infty)\times\RR^N\,.
\eeqn
Equivalently, $v(\tau,y)=u\left( e^{q\tau}-1,y e^\tau \right)$ for $(\tau,y)\in [0,\infty)\times\RR^N$ and it follows from \eqref{a1}, \eqref{a2} that $v$ solves
\bear
\label{d2}
\partial_\tau v & = & y\cdot\nabla v + q\ \vert\nabla v\vert^q + q\ e^{-(p-q)\tau}\ \Delta_p v \,, \quad (\tau,y)\in (0,\infty)\times\RR^N\,,\\
\label{d3}
v(0) & = & u_0\,, \quad y\in\RR^N\,.
\eear
We also infer from (\ref{b1}) and (\ref{b2}) that there is a positive constant $C_5(u_0)$ depending only on $N$, $p$, $q$, and $u_0$ such that
\beqn
\label{d4}
\Vert v(\tau)\Vert_\infty + \Vert\nabla v(\tau)\Vert_\infty  \le C_5(u_0) \,,\qquad \tau\ge 0\,,
\eeqn 
while \eqref{a4} reads 
\beqn
\label{d6}
\lim_{\tau\to\infty} \|v(\tau)\|_\infty=M_\infty>0\,.
\eeqn

Formally, since $p>q$, the diffusion term vanishes in the large time limit and we expect the large time behaviour of the solution $v$ to \eqref{d2}, \eqref{d3} to look like that of the solutions to the first-order Hamilton-Jacobi equation 
\beqn
\label{d6b}
\partial_\tau w - y\cdot\nabla w - q\ \vert\nabla w\vert^q=0 \;\;\mbox{ in }\;\; Q_\infty\,.
\eeqn
Now, to investigate the large time behaviour of first-order Hamilton-Jacobi equations, an efficient approach has been developed in \cite{NR99,Ro01} which relies on the relaxed half-limits method introduced in \cite{BP88}. More precisely, for $(\tau,y)\in (0,\infty)\times\RR^N$, we define the relaxed half-limits $v_*$ and $v^*$ by
\beqn 
\label{d7} 
v_*(y):= \liminf_{(\sigma,z,\lambda)\to(\tau,y,\infty)}{v(\sigma+\lambda,z)} 
\;\;\mbox{ and }\;\; v^*(y):= 
\limsup_{(\sigma,z,\lambda)\to(\tau,y,\infty)}{v(\sigma+\lambda,z)} \,.
\eeqn 
These relaxed half-limits are well-defined thanks to \eqref{d4} and we first note that the right-hand sides of the above definitions 
indeed do not depend on $\tau>0$. In addition, 
\beqn 
\label{d8} 
0 \le v_*(x) \le v^*(x) \le M_\infty \;\;\mbox{ for }\;\; y\in\RR^N 
\eeqn 
by \eqref{d6}, while \eqref{d4} and the Rademacher theorem ensure that $v_*$ and $v^*$ both belong to $W^{1,\infty}(\RR^N)$. Finally, by \cite[Th\'eor\`eme~4.1]{Bl94}  applied to equation \eqref{d2}, $v^*$ and $v_*$ are viscosity subsolution and supersolution, respectively, to the Hamilton-Jacobi equation 
\beqn 
\label{d9} 
\mathcal{H}(y,\nabla w) := -y\cdot\nabla w - q\ |\nabla w|^q = 0 \;\;\;\mbox{ in }\;\;\; \RR^N\,. 
\eeqn 

We now aim at showing that $v^*$ and $v_*$ coincide. However, the equation \eqref{d9} has infinitely many solutions as $y\longmapsto  \left( c - \gamma_q\ \vert y\vert^{q/(q-1)} \right)_+$ solves \eqref{d9} for any $c>0$ .The information obtained so far on $v_*$ and $v^*$ are thus not sufficient and are supplemented by the next two results.

\begin{lem}\label{led1}
Given $\e\in (0,1)$, there is $R_\e>1/\e$ such that
\beqn
\label{d10}
v(\tau,y) \le \e \;\;\mbox{ for }\;\; \tau\ge 0 \;\;\mbox{ and }\;\; y\in \RR^N\setminus B(0,R_\e)\,,
\eeqn
and $0\le v_*(y)\le v^*(y)\le \e$ for $y\in \RR^N\setminus B(0,R_\e)$.
\end{lem}

In other words, $v(\tau)$ belongs to $\mathcal{C}_0(\RR^N)$ for each $\tau\ge 0$ in a way which is uniform with respect to $\tau\ge 0$. 

\noindent\textbf{Proof.} We first construct a supersolution to (\ref{d2}) in $(0,\infty)\times\RR^N\setminus B(0,R)$ for $R$ large enough. To this end, consider $R\ge R_c:= 1 + \left( q\ (2\ \|u_0\|_\infty)^{q-1} + 3pq\ (2\ \|u_0\|_\infty)^{p-2} \right)^{1/q}$ and put $\Sigma_R(y) = \Vert u_0\Vert_\infty\ R^2\ \vert y\vert^{-2}$ for $y\in\RR^N\setminus B(0,R)$. Let $\mathcal{L}$ be the parabolic operator defined by
$$
\mathcal{L} w(\tau,y) := \partial_\tau w(\tau,y) - y\cdot\nabla w(\tau,y) - q\ \vert\nabla w(\tau,y)\vert^q - q\ e^{-(p-q)\tau}\ \Delta_p w(\tau,y)
$$
for $(\tau,y)\in Q_\infty$ (so that $\mathcal{L}v=0$ by \eqref{d2}). Then, if  $y\in \RR^N\setminus B(0,R)$, we have
\bean
\mathcal{L} \Sigma_R(y) & = & 2\ \Sigma_R(y) - q\ \frac{2^q}{|y|^q}\  \Sigma_R(y)^q + q\ 2^{p-1}\ \frac{N+2-3p}{|y|^p}\ \Sigma_R(y)^{p-1}\ e^{-(p-q)\tau} \\
& \ge & 2\ \Sigma_R(y)\ \left\{ 1 - q\ (2\ \Vert u_0\Vert_\infty)^{q-1}\ \frac{R^{2(q-1)}}{\vert y\vert^{3q-2}} - 3pq\  e^{-(p-q)\tau}\ (2\ \Vert u_0\Vert_\infty)^{p-2}\ \frac{R^{2(p-2)}}{\vert y\vert^{3p-4}}  \right\} \\
& \ge & 2\ \Sigma_R(y)\ \left\{ 1 - q\ (2\ \Vert u_0\Vert_\infty)^{q-1}\ R^{-q} - 3pq\  e^{-(p-q)\tau}\ (2\ \Vert u_0\Vert_\infty)^{p-2}\ R^{-p} \right\} \\
& \ge & 0
\eean
by the choice of $R$. Consequently, $\Sigma_R$ is a supersolution to (\ref{d2}) in $(0,\infty)\times\RR^N\setminus B(0,R)$ for $R\ge R_c$.

Now, fix $\e\in (0,1)$. Since $u_0\in\mathcal{C}_0(\RR^N)$, there is $\rho_\e\ge  \max{\{ R_c , \e^{-1} \}}$ such that $u_0(y)\le \e/2$ if $\vert y\vert\ge \rho_\e$. We then infer from the monotonicity of $\Sigma_R$ and \eqref{b1} that
$$
u_0(y) - \frac{\e}{2} - \Sigma_{\rho_\e}(z) \le - \Sigma_{\rho_\e}(z) \le 0
$$
if $|y|\ge\rho_\e$ and $|z|\ge\rho_\e$, 
$$
v(\tau,y) - \frac{\e}{2} - \Sigma_{\rho_\e}(z) \le \|u_0\|_\infty - \Sigma_{\rho_\e}(z) = \Sigma_{\rho_\e}(y) - \Sigma_{\rho_\e}(z) \le \frac{2\ \|u_0\|_\infty}{\rho_\e}\ |y-z|
$$
if $|y|=\rho_\e$, $|z|\ge\rho_\e$, and $\tau\ge 0$, and
$$
v(\tau,y) - \frac{\e}{2} - \Sigma_{\rho_\e}(z) \le \|u_0\|_\infty - \|u_0\|_\infty \le 0
$$
if $|y|\ge\rho_\e$, $|z|=\rho_\e$, and $\tau\ge 0$. As $v-\e/2$ and $\Sigma_{\rho_\e}$ are subsolution and supersolution, respectively, to \eqref{d2}, the comparison principle \cite[Theorem~4.1]{GGIS91} warrants that $v(\tau,y) - \e/2 \le \Sigma_{\rho_\e}(y)$ for $\tau\ge 0$ and $|y|\ge\rho_\e$. It remains to choose $R_\e\ge\rho_\e$ such that $\Sigma_{\rho_\e}(y)\le \e/2$ for $|y|\ge R_\e$ to complete the proof of \eqref{d10}. The last assertion of Lemma~\ref{led1} is then a straightforward consequence of the definition \eqref{d7} and \eqref{d10}. \qed

\medskip

We next use the semiconvexity estimate \eqref{a11} (and more precisely its consequence stated in Corollary~\ref{coc4}) to show that $v_*$ lies above the profile $H_\infty$ defined in \eqref{a6}.

\begin{lem}\label{led2}
For $y\in\RR^N$, we have
\beqn
\label{d11}
H_\infty(y) \le v_*(y) \le v^*(y)\,.
\eeqn
\end{lem}

\noindent\textbf{Proof.} For $\tau\ge 0$, $y\in\RR^N$, $\xi_0\in\RR$ and $\xi\in\RR^N$, we set $\mathcal{F}(\tau,y,\xi_0,\xi) := \xi_0 - y\cdot\xi - q\ \vert\xi\vert^q + \kappa_2\ e^{-(p-q)\tau}$ with $\kappa_2:= q\ \kappa_1\ e^q/(e^q-1)$, the constant $\kappa_1$ being defined in \eqref{a11}. It then readily follows from Corollary~\ref{coc4} that 
\beqn
\label{d12}
v \;\;\mbox{ is a supersolution to }\;\; \mathcal{F}(\tau,y,\partial_\tau w,\nabla w)=0 \;\;\mbox{ in }\;\;  (1,\infty)\times\RR^N\,.
\eeqn

We next fix $\tau_0>1$ and denote by $V$ the (viscosity) solution to 
\bean
\partial_\tau V - y\cdot \nabla V - q\ \vert\nabla V\vert^q & = & 0 \,, \qquad (\tau,y)\in (\tau_0,\infty)\times\RR^N\,,\\
V(\tau_0) & = & v(\tau_0)\,, \qquad y\in \RR^N\,.
\eean
On the one hand, a straightforward computation shows that the function $\tilde{V}$ defined by
$$
\tilde{V}(\tau,y) := V(\tau,y) - \kappa_2\ \int_{\tau_0}^\tau e^{-(p-q)s}\ ds\,, \qquad (\tau,y)\in (\tau_0,\infty)\times\RR^N\,,
$$
is the (viscosity) solution to $\mathcal{F}(\tau,y,\partial_\tau \tilde{V},\nabla \tilde{V})=0$ in $(\tau_0,\infty)\times\RR^N$ with initial condition $\tilde{V}(\tau_0)=v(\tau_0)$. Recalling (\ref{d12}), we infer from the comparison principle that 
\beqn
\label{d13}
\tilde{V}(\tau,y)\le v(\tau,y) \;\;\mbox{ for }\;\; (\tau,y)\in (\tau_0,\infty)\times\RR^N\,.
\eeqn 
On the other hand, it follows from Proposition~\ref{prap1} that
$$
\lim_{\tau\to\infty} \sup_{y\in\RR^N}{ \left\vert V(\tau,y) - \left( \Vert v(\tau_0)\Vert_\infty - \gamma_q\  \vert y\vert^{q/(q-1)} \right)_+ \right\vert } = 0\,.
$$
We may then pass to the limit as $\tau\to\infty$ in (\ref{d13}) and use the definition \eqref{d7} to conclude that
$$
\left( \Vert v(\tau_0)\Vert_\infty - \gamma_q\  \vert y\vert^{q/(q-1)} \right)_+ - \kappa_2\ \int_{\tau_0}^\infty e^{-(p-q)s}\ ds \le v_*(y)\le v^*(y)
$$
for $y\in\RR^N$. Letting $\tau_0\to\infty$ in the above inequality with the help of \eqref{d6} completes the proof of the lemma. \qed

\medskip

We are now in a position to complete the proof of Theorem~\ref{tha1}. To this end, fix $\varepsilon\in (0,1)$. Lemma~\ref{led1} ensures that $v^*(y)\le\e$ for $\vert y\vert\ge R_\e\ge 1/\e$ while the continuity of $H_\infty$ implies that there is $r_\e\in (0,\e)$ such that $H_\infty(y)\ge M_\infty - \e$ for $|y|\le r_\e$. Recalling \eqref{d8}, we realize that
\beqn
\label{d14}
\left\{
\begin{array}{lcl}
v^*(y)-\e \le 0 \le H_\infty(y) & \mbox{ if } & \vert y\vert=R_\e \,,\\
& & \\
v^*(y)-\e \le M_\infty-\e \le H_\infty(y) & \mbox{ if } & \vert y\vert=r_\e\,.
\end{array}
\right.
\eeqn
Moreover, introducing $\psi(y)=-\gamma_q\ \vert y\vert^{q/(q-1)}/2$, we have
\beqn
\label{d15}
\mathcal{H}(y,\nabla\psi(y)) = \frac{q \gamma_q}{2(q-1)}\ \vert y\vert^{q/(q-1)}\ \left( 1 - \frac{1}{2^{q-1}} \right) > 0 \;\;\mbox{ if }\;\; r_\e < \vert y\vert < R_\e\,,
\eeqn
the Hamiltonian $\mathcal{H}$ being defined in \eqref{d9}.
Summarizing, we have shown that $H_\infty$ and $v^*-\e$ are supersolution and subsolution, respectively, to \eqref{d9} in $\Omega_\e:=\left\{ y\in\RR^N\; : \; r_\e < \vert y\vert < R_\e\right\}$ with $v^*-\e\le H_\infty$ on $\partial\Omega_\e$ by \eqref{d14}. Owing to \eqref{d15} and the concavity of $\mathcal{H}$ with respect to its second variable, we may apply \cite[Theorem~1]{I87} to conclude that $v^*-\e\le H_\infty$ in $\Omega_\e$. This property being valid for each $\e\in (0,1)$, we actually have $v^*\le H_\infty$ in $\RR^N$ by passing to the limit as $\e\to 0$ thanks to the properties of $r_\e$ and $R_\e$. Recalling (\ref{d11}), we have thus established that $v^*=v_*=H_\infty$ in $\RR^N$. In particular, the property $v^*=v_*$ and the definition \eqref{d7} provide the uniform convergence of $\{ v(\tau)\}_{\tau\ge 0}$ towards $v^*=H_\infty$ on every compact subset of $\RR^N$ as $\tau\to\infty$, see \cite[Lemme~4.1]{Bl94} or \cite[Lemma~V.1.9]{BdCD97}. Combining this local convergence with Lemma~\ref{led1} actually gives
\beqn
\label{d16} 
\lim_{\tau\to\infty} \|v(\tau)-H_\infty\|_\infty=0\,. 
\eeqn 
Theorem~\ref{tha1} then readily follows after writing the convergence (\ref{d16}) in the original variables $(t,x)$ for the function $u$ and noticing that $\Vert h_\infty(1+t) - h_\infty(t) \Vert_\infty\longrightarrow 0$ as $t\to\infty$. \qed

\section{Limit value of $\|u(t)\|_\infty$}\label{lv}

\setcounter{thm}{0}
\setcounter{equation}{0}

This section is devoted to the proof of Proposition~\ref{pra3}, for which three cases are to be distinguished and handled differently: $q\in (1,p-1]$, $q\in (p-1,q_\star]$, and $q\in (q_\star,p)$.

\medskip

\noindent\textbf{Proof of Proposition~\ref{pra3}: $q\in (1,p-1]$.} We proceed as in \cite[Proposition~1]{LS05} (where a similar result is proved for $p=2$ and $q=1$). For $\alpha>N/2$, $\delta>0$, and $x\in\RR^N$, we set $\varrho_\delta(x):= \left( 1 + \delta\ |x|^2 \right)^{-\alpha}$. Clearly, $\varrho_\delta\in L^1(\RR^N)$ and it follows from \eqref{b3} that
\bean
\frac{d}{dt} \int_{\RR^N} \varrho_\delta(x)\ u(t,x)\ dx & = & \int_{\RR^N} \left\{ \varrho_\delta(x)\ |\nabla u(t,x)|^q - |\nabla u(t,x)|^{p-2}\ \nabla u(t,x) \cdot \nabla \varrho_\delta(x) \right\}\ dx \\
& \ge & \int_{\RR^N} \varrho_\delta(x)\ |\nabla u(t,x)|^q\ \left( 1 - |\nabla u(t,x)|^{p-1-q}\ \frac{|\nabla \varrho_\delta(x)|}{\varrho_\delta(x)} \right)\ dx\,.
\eean
Recalling that $\|\nabla u(t)\|_\infty\le \|\nabla u_0\|_\infty$  by \eqref{a3} and \eqref{b2} and noticing that $|\nabla\varrho_\delta|\le\alpha\ \delta^{1/2}\ \varrho_\delta$, we further obtain
$$
\frac{d}{dt} \int_{\RR^N} \varrho_\delta(x)\ u(t,x)\ dx \ge \int_{\RR^N} \varrho_\delta(x)\ |\nabla u(t,x)|^q\ \left( 1 - \alpha\ \delta^{1/2}\ \|\nabla u_0\|_\infty^{p-1-q} \right)\ dx\,.
$$
Choosing $\delta = \|\nabla u_0\|_\infty^{2(q+1-p)}/\alpha^2$ and integrating with respect to time give
$$
\|u(t)\|_\infty\ \|\varrho_\delta\|_1 \ge \int_{\RR^N} \varrho_\delta(x)\ u(t,x)\ dx \ge \int_{\RR^N} \varrho_\delta(x)\ u_0(x)\ dx > 0\,.
$$
We then pass to the limit as $t\to\infty$ to conclude that $M_\infty>0$. \qed

\medskip

We next turn to the case $q\in (p-1,q_\star]$ which turns out to be more complicated and requires two preparatory results. 

\begin{lem}\label{lee1}
Assume that $q\in (1,q_\star]$ and let $u$ be a non-negative solution to \eqref{a1}, \eqref{a2} with a compactly supported initial condition $u_0$ satisfying \eqref{a3}.  Then $u(t)\in L^1(\RR^N)$ for each $t\ge 0$, the function $t\longmapsto \|u(t)\|_1$ is non-decreasing and
\beqn
\label{e1}
\lim_{t\to\infty} \|u(t)\|_1 = \infty\,.
\eeqn
\end{lem}

\noindent\textbf{Proof.} For every $t\ge 0$, $u(t)$ is bounded and compactly supported by \eqref{b1} and Proposition~\ref{prb2}, and is thus in $L^1(\RR^N)$. The time monotonicity of the $L^1$-norm of $u$ then readily follows from \eqref{b3} with $\vartheta=1$, a valid choice in this particular case as $u(t)$ is compactly supported. It further follows from \eqref{b3} with $\vartheta=1$ that
\beqn
\label{e2}
\|u(t)\|_1 \ge \|u(T)\|_1 + \int_T^t \|\nabla u(s)\|_q^q\ ds \;\;\mbox{ for } t>T\ge 0\,.
\eeqn

Consider next $T>0$ and $t>T$. Recalling the Gagliardo-Nirenberg inequality
\beqn
\label{e3}
\|w\|_q \le C_6\ \|\nabla w\|_q^{N(q-1)/(N(q-1)+q)}\ \|w\|_1^{q/(N(q-1)+q)}\,, \quad w\in W^{1,q}(\RR^N)\cap L^1(\RR^N)\,,
\eeqn
we infer from \eqref{e2}, \eqref{e3}, and the time monotonicity of the $L^1$-norm of $u$ that 
\bean
\|u(t)\|_1^{1 +(q^2/N(q-1))} & \ge & \|u(t)\|_1^{q^2/N(q-1)}\ \left( \|u(T)\|_1 + \int_T^t \|\nabla u(s)\|_q^q\ ds \right) \\
& \ge & \int_T^t \|u(s)\|_1^{q^2/N(q-1)}\ \|\nabla u(s)\|_q^q\ ds \\
& \ge & C_7\ \int_T^t \left( \|u(s)\|_q^q \right)^{(N(q-1)+q)/N(q-1)}\ ds\,.
\eean
If $\varphi$ denotes the solution to the $p$-Laplacian equation $\partial_t \varphi - \Delta_p \varphi = 0$ in $Q_\infty$ with initial condition $\varphi(0)=u_0$, the comparison principle readily implies that 
\beqn
\label{e5}
\varphi(t,x)\le u(t,x)\,, \quad (t,x)\in Q_\infty\,.
\eeqn
Inserting this estimate in the previous lower bound for $\|u(t)\|_1$, we end up with
\beqn
\label{e4}
\|u(t)\|_1^{1 +(q^2/N(q-1))} \ge C_7\ \int_T^t \left( \| \varphi(s)\|_q^q \right)^{(N(q-1)+q)/N(q-1)}\ ds\,.
\eeqn
Now, by Proposition~\ref{prb3} we have
$$
\lim_{s\to\infty} s^{N(q-1)/(N(p-2)+p)}\ \left\| \varphi(s) - \mathcal{B}_{\|u_0\|_1}(s) \right\|_q^q = 0
$$
and
$$
\left\|\mathcal{B}_{\|u_0\|_1}(s) \right\|_q^q = C_8\ s^{-N(q-1)/(N(p-2)+p)}\,,
$$
so that
\bean
\| \varphi(s)\|_q^q & \ge & \left( \left\|\mathcal{B}_{\|u_0\|_1}(s) \right\|_q - \left\| \varphi(s) - \mathcal{B}_{\|u_0\|_1}(s) \right\|_q \right)^q\\
& \ge & s^{-N(q-1)/(N(p-2)+p)}\ \left( C_8 - \left\| \varphi(s) - \mathcal{B}_{\|u_0\|_1}(s) \right\|_q \right)^q\\
& \ge & \left( \frac{C_8}{2} \right)^q\ s^{-N(q-1)/(N(p-2)+p)}
\eean
for $s\ge T$, provided $T$ is chosen sufficiently large. Inserting this estimate in \eqref{e4} gives
\bean
\|u(t)\|_1^{1 +(q^2/N(q-1))} & \ge & C_9\ \int_T^t s^{-N(q-1)/(N(p-2)+p)}\ ds \\
& \ge & C_{10}\ \left\{
\begin{array}{lcl}
\displaystyle{t^{(N+1)(q_\star-q)/(N(p-2)+p)} - T^{(N+1)(q_\star-q)/(N(p-2)+p)}} & \mbox{ if } & q\in (1,q_\star)\,,\\
 & & \\
\log(t/T) & \mbox{ if } & q=q_\star\,.
\end{array}
\right.
\eean
We then let $t\to\infty$ to obtain the claimed result. \qed

\medskip

We next argue as in \cite[Lemma~14]{Gi05} (for $p=2$) to show that, if $q\in (p-1,p)$ and $M_\infty=0$, then the $L^\infty$-norm of $u(t)$ decays faster than an explicit rate.

\begin{lem}\label{lee2}
Assume that $q\in (p-1,p)$ and let $u$ be a non-negative solution to \eqref{a1}, \eqref{a2} with an initial condition $u_0$ satisfying \eqref{a3}.  If $M_\infty=0$ in \eqref{a4}, then 
\beqn
\label{e6}
\|u(t)\|_\infty \le C_{11}\ t^{-(p-q)/(2q-p)} \;\;\mbox{ for }\;\; t>0\,.
\eeqn
\end{lem}

Observe that the assumptions $p>2$ and $q\in (p-1,p)$ imply that $2q>p$ and $(p-q)/(2q-p)>0$. 

\noindent\textbf{Proof.} Consider a non-negative function $\eta\in \mathcal{C}^\infty(\RR^N)$ with compact support in $B(0,1)$ and $\| \eta\|_1=1$. We then define a sequence of mollifiers $(\eta_\delta)_\delta$ by $\eta_\delta(x) := \eta(x/\delta)/\delta^N$ for $x\in\RR^N$ and $\delta\in (0,1)$. For $(t,x_0)\in Q_\infty$ and $T>t$, we take $\vartheta(x)= \eta_\delta(x-x_0)$ in \eqref{b3} and infer from \eqref{a11} (with $s=t/2$) that
\bean
\|u(T)\|_\infty & \ge & \int_{\RR^N} u(T,x)\ \eta_\delta(x-x_0)\ dx \\
& \ge & \int_{\RR^N} u(t,x)\ \eta_\delta(x-x_0)\ dx - \int_t^T \int_{\RR^N} |\nabla u(s,x)|^{p-2}\ \nabla u(s,x) \cdot \nabla \eta_\delta(t,x-x_0)\ dxds \\
& \ge & \int_{\RR^N} u(t,x)\ \eta_\delta(x-x_0)\ dx - 2^{p/q}\ \kappa_1\ \left\| u\left( \frac{t}{2} \right) \right\|_\infty^{(p-q)/q}\ \int_t^T (2s-t)^{-p/q}\ ds \\
& \ge & \int_{\RR^N} u(t,x)\ \eta_\delta(x-x_0)\ dx - C_{12}\ \left\| u\left( \frac{t}{2} \right) \right\|_\infty^{(p-q)/q}\ \left( t^{(q-p)/p} - T^{(q-p)/p} \right)\,.
\eean
Owing to the continuity of $u$, we may pass to the limit as $\delta\to 0$ in the above inequality and deduce that 
$$
\|u(T)\|_\infty \ge u(t,x_0) - C_{12}\ \left\| u\left( \frac{t}{2} \right) \right\|_\infty^{(p-q)/q}\ \left( t^{(q-p)/p} - T^{(q-p)/p} \right)\,.
$$
But the above inequality is valid for all $x_0\in\RR^N$ and we thus end up with
$$
\|u(T)\|_\infty \ge \|u(t)\|_\infty - C_{12}\ \left\| u\left( \frac{t}{2} \right) \right\|_\infty^{(p-q)/q}\ \left( t^{(q-p)/p} - T^{(q-p)/p} \right)\,.
$$
Finally, as $q<p$, we may let $T\to\infty$ in the previous inequality and use the assumption $M_\infty=0$ to conclude that
$$
\|u(t)\|_\infty \le C_{12}\ \left\| u\left( \frac{t}{2} \right) \right\|_\infty^{(p-q)/q}\ t^{(q-p)/p}\,,
$$
or, equivalently, as $2q>p$,
$$
t^{(p-q)/(2q-p)}\ \|u(t)\|_\infty \le C_{13}\ \left\{ \left( \frac{t}{2} \right)^{(p-q)/(2q-p)}\ \left\| u\left( \frac{t}{2} \right) \right\|_\infty\right\}^{(p-q)/q}
$$
for $t\ge 0$. Introducing
$$
A(t):=\sup_{s\in (0,t)}{ \left\{ s^{(p-q)/(2q-p)}\ \| u(s)\|_\infty \right\}} \in [0,\infty)\,, \quad t\ge 0\,,
$$
we deduce from the previous inequality that $A(t)\le C_{13}\ A(t)^{(p-q)/q}$, whence $A(t)\le C_{13}^{q/(2q-p)}$ for $t\ge 0$. This bound being valid for each $t>0$, the proof of \eqref{e6} is complete. \qed

\medskip

\noindent\textbf{Proof of Proposition~\ref{pra3}: $q\in (p-1,q_\star]$.}

\smallskip

\noindent\textsl{Step~1:} We first consider a compactly supported initial condition $u_0$ satisfying \eqref{a3} and assume for contradiction that $M_\infty=0$. On the one hand, according to Lemma~\ref{lee2} and the assumption $q\le q_\star$, there holds
\beqn
\label{e7}
\limsup_{t\to\infty} t^{N/(N(p-2)+p)}\ \|u(t)\|_\infty  \le C_{11}\ t^{(N(p+1)(q-q_\star))/((2q-p)(N(p-2)+p))} \le C_{11}\,.
\eeqn
On the other hand, fix $t_0>0$ and let $\varphi$ be the solution to the $p$-Laplacian equation $\partial_t \varphi - \Delta_p \varphi = 0$ in $Q_\infty$ with initial condition $\varphi(0)=u(t_0)$. As $u_0$ is compactly supported, so is $u(t_0)$ by Proposition~\ref{prb2} and $u(t_0)$ thus belongs to $L^1(\RR^N)$. Moreover, the comparison principle warrants that $u(t,x)\ge \varphi(t-t_0,x)$ for $(t,x)\in [t_0,\infty)\times\RR^N$. We then infer from the above properties and Proposition~\ref{prb3} that, for $t>t_0$, 
\bean
t^{N/(N(p-2)+p)}\ \|u(t)\|_\infty  & \ge & (t-t_0)^{N/(N(p-2)+p}\ \|u(t)\|_\infty  \\
& \ge & (t-t_0)^{N/(N(p-2)+p)}\ \| \varphi(t-t_0)\|_\infty  \\
& \ge & (t-t_0)^{N/(N(p-2)+p)}\ \left\| \mathcal{B}_{\|u(t_0)\|_1}(t-t_0) \right\|_\infty  \\
& & -\ (t-t_0)^{N/(N(p-2)+p)}\ \left\| \mathcal{B}_{\|u(t_0)\|_1}(t-t_0) - \varphi(t-t_0) \right\|_\infty  \\
& \ge & C_{14}\ \|u(t_0)\|_1^{N/(N(p-2)+p)}   \\
& & -\ (t-t_0)^{N/(N(p-2)+p)}\ \left\| \mathcal{B}_{\|u(t_0)\|_1}(t-t_0) - \varphi(t-t_0) \right\|_\infty\,.
\eean
Using once more Proposition~\ref{prb3}, we may pass to the limit as $t\to\infty$ in the previous inequality to obtain
\beqn
\label{e8}
\liminf_{t\to\infty} t^{N/(N(p-2)+p)}\ \|u(t)\|_\infty  \ge C_{14}\ \|u(t_0)\|_1^{N/(N(p-2)+p)}\,.
\eeqn
Combining \eqref{e7} and \eqref{e8} yields $\|u(t_0)\|_1\le C_{15}$ for all $t_0>0$ which contradicts Lemma~\ref{lee1}. Therefore, $M_\infty>0$. 

\smallskip

\noindent\textsl{Step~2:} Now, if $u_0$ is an arbitrary initial condition satisfying \eqref{a3}, there clearly exists a compactly supported initial condition $\tilde{u}_0$ satisfying \eqref{a3} and such that $u_0\ge \tilde{u}_0$ in $\RR^N$. Introducing the solution $\tilde{u}$ to \eqref{a1} with initial condition $\tilde{u}_0$, the comparison principle entails that $u\ge\tilde{u}$ in $Q_\infty$, hence
$$
M_\infty\ge \lim_{t\to\infty} \|\tilde{u}(t)\|_\infty\,.
$$
The first step of the proof ensures that the right-hand side of the above inequality is positive which completes the proof. \qed

\medskip

It remains to investigate the case $q\in (q_\star,p)$, for which we adapt the proof of \cite[Theorem~2.4(b)]{BKL04}.

\noindent\textbf{Proof of Proposition~\ref{pra3}: $q\in (q_\star,p)$.} We put
$$
m_0:=\left| \inf_{y\in\RR^N} {\Delta_p u_0(y)} \right|\,.
$$
As in the proof of Lemma~\ref{lee1}, let $\eta\in \mathcal{C}^\infty(\RR^N)$ be a non-negative function with compact support in $B(0,1)$ and $\|\eta\|_1=1$, and define a sequence of mollifiers $(\eta_\delta)_\delta$ by $\eta_\delta(x) := \eta(x/\delta)/\delta^N$ for $x\in\RR^N$ and $\delta\in (0,1)$. For $(t,x_0)\in Q_\infty$ and $T\in (0,t)$, we take $\vartheta(x)= \eta_\delta(x-x_0)$ in \eqref{b3} and infer from \eqref{a11} (with $s=0$) and Corollary~\ref{coc3} that
\bean
\|u(t)\|_\infty & \ge & \int_{\RR^N} u(t,x)\ \eta_\delta(x-x_0)\ dx \\
& \ge & \int_{\RR^N} u_0(x)\ \eta_\delta(x-x_0)\ dx - \int_0^t \int_{\RR^N} |\nabla u(s,x)|^{p-2}\ \nabla u(s,x) \cdot \nabla \eta_\delta(t,x-x_0)\ dxds \\
& \ge & \int_{\RR^N} u_0(x)\ \eta_\delta(x-x_0)\ dx - \int_0^T m_0\ ds - \kappa_1\ \|u_0\|_\infty^{(p-q)/q}\ \int_T^t s^{-p/q}\ ds \\
& \ge & \int_{\RR^N} u_0(x)\ \eta_\delta(x-x_0)\ dx - T\ m_0 - C_{16}\ \|u_0\|_\infty^{(p-q)/q}\ \left( T^{(q-p)/p} - t^{(q-p)/p} \right)\,.
\eean
Owing to the continuity of $u_0$, we may pass to the limit as $\delta\to 0$ in the above inequality and deduce that 
$$
\|u(t)\|_\infty \ge u_0(x_0) - T\ m_0 - C_{16}\ \|u_0\|_\infty^{(p-q)/q}\ \left( T^{(q-p)/p} - t^{(q-p)/p} \right)\,.
$$
Since $q<p$, we may let $t\to\infty$ in the above inequality and take the supremum with respect to $x_0$ to conclude that
$$
M_\infty \ge \|u_0\|_\infty - T\ m_0 - C_{16}\ \|u_0\|_\infty^{(p-q)/q}\ T^{(q-p)/p}\,.
$$
Next, for $\beta\in (0,1)$, the choice $T=\|u_0\|_\infty^{(p-q)/q}\ (\beta+m_0)^{-q/p}$ in the previous inequality yields
$$
M_\infty \ge \|u_0\|_\infty^{(p-q)/q}\ \left( \|u_0\|_\infty^{p/q} - (1+C_{16})\ (\beta+m_0)^{(p-q)/p} \right)\,.
$$
This inequality being valid for every $\beta\in (0,1)$, we conclude that 
$$
M_\infty \ge \|u_0\|_\infty^{(p-q)/q}\ \left( \|u_0\|_\infty^{p/q} - (1+C_{16})\ m_0^{(p-q)/p} \right)>0
$$
as soon as \eqref{a10} is fulfilled with $\kappa_0=(1+C_{16})^{p/q}$.
\qed

\appendix
\section{Convergence for the Hamilton-Jacobi equation \eqref{d6b}}

\setcounter{thm}{0}
\setcounter{equation}{0}

In this section, we study the large behaviour of non-negative solutions to the Hamilton-Jacobi equation \eqref{d6b} with initial data in $\mathcal{C}_0(\RR^N)$ and show their convergence to a steady state uniquely determined by the $L^\infty$-norm of the initial data. Though the large time behaviour of solutions to first-order Hamilton-Jacobi equations has received considerable attention in recent years (see \cite{BR06,BS00,I08,NR99,Ro01} and the references therein), the particular case of \eqref{d6b} does not seem to have been investigated in the literature. We thus provide a simple proof relying on the Hopf-Lax-Oleinik formula.

\begin{prop}\label{prap1}
Let $q>1$ and consider a non-negative function $h_0\in \mathcal{C}_0(\RR^N)$. Let $h$ be the unique viscosity solution to the Cauchy problem
\bear
\label{ap1}
\partial_\tau h - y\cdot \nabla h - q\ \vert\nabla h\vert^q & = & 0 \,, \qquad (\tau,y)\in (0,\infty)\times\RR^N\,,\\
\label{ap2}
h(0) & = & h_0\,, \qquad y\in \RR^N\,.
\eear
Then
\beqn
\label{ap3}
\lim_{\tau\to\infty} \left\| h(\tau) - h_s \right\|_\infty = 0
\eeqn
with
$$
h_s(y) := \left( \|h_0\|_\infty - \gamma_q\  \vert y\vert^{q/(q-1)} \right)_+\,, \qquad y\in\RR^N\,,
$$
the constant $\gamma_q=(q-1)\ q^{-q/(q-1)}$ being defined in Theorem~\ref{tha1}.
\end{prop}

Thanks to the concavity of the Hamiltonian $\mathcal{H}(y,\xi) = -y\cdot\xi - q\ |\xi|^q$, $(y,\xi)\in\RR^N\times\RR^N$, with respect to its second variable, the Hopf-Lax-Oleinik formula provides a representation formula  for the solution $h$ to \eqref{ap1}, \eqref{ap2} which can be used to prove \eqref{ap3}. 

\noindent\textbf{Proof.} We first recall that $h$ is given by the Hopf-Lax-Oleinik formula 
$$
h(\tau,y) = \sup_{z\in\RR^N}{\left\{ h_0(z) - \gamma_q\ \left| y-z\ e^{-\tau} \right|^{q/(q-1)}\ \left( 1 - e^{-q\tau} \right)^{-1/(q-1)} \right\}} 
$$
for $(\tau,y)\in [0,\infty)\times\RR^N$, see, e.g., \cite[Chapter~3]{Ev98}. Since $h(\tau,y)\ge h_0(y e^\tau)\ge 0$, we have in fact
$$
h(\tau,y) = \sup_{z\in\RR^N}{\left\{ \left( h_0(z) - \gamma_q\ \left| y-z\ e^{-\tau} \right|^{q/(q-1)}\ \left( 1 - e^{-q\tau} \right)^{-1/(q-1)} \right)_+ \right\}} 
$$
for $(\tau,y)\in [0,\infty)\times\RR^N$. 

Consider now $\beta\in (0,1)$. As $h_0\in\mathcal{C}_0(\RR^N)$, there is $R_\beta>(\Vert h_0\Vert_\infty/\gamma_q)^{(q-1)/q}$ such that 
\beqn
\label{ap4}
h_0(z) \le \beta \quad\mbox{ for }\quad |z|\ge R_\beta\,.
\eeqn 
On the one hand, if $(\tau,y)\in [\log{R_\beta},\infty)\times\RR^N$ and $z\in\RR^N$, we have either $|z|\ge R_\beta$ and 
\bean
& & \left| \left( h_0(z) - \gamma_q\ \left| y-z\ e^{-\tau} \right|^{q/(q-1)}\ \left( 1 - e^{-q\tau} \right)^{-1/(q-1)} \right)_+ - \left( h_0(z) - \gamma_q\ \left| y \right|^{q/(q-1)}\  \right)_+\right| \\
& \le & \left( h_0(z) - \gamma_q\ \left| y-z\ e^{-\tau} \right|^{q/(q-1)}\ \left( 1 - e^{-q\tau} \right)^{-1/(q-1)} \right)_+ + \left( h_0(z) - \gamma_q\ \left| y \right|^{q/(q-1)} \right)_+ \\
& \le & 2\ \beta
\eean
by \eqref{ap4} or $z\in B(0,R_\beta)$ and 
\bean
& & \left| \left( h_0(z) - \gamma_q\ \left| y-z\ e^{-\tau} \right|^{q/(q-1)}\ \left( 1 - e^{-q\tau} \right)^{-1/(q-1)} \right)_+ - \left( h_0(z) - \gamma_q\ \left| y \right|^{q/(q-1)}\  \right)_+\right| \\
& \le & \gamma_q\ \left| y-z\ e^{-\tau} \right|^{q/(q-1)}\ \left\{ \left( 1 - e^{-q\tau} \right)^{-1/(q-1)} - 1 \right\} + \gamma_q\ \left| \left| y-z\ e^{-\tau} \right|^{q/(q-1)} - \left| y \right|^{q/(q-1)} \right| \\
& \le &  \gamma_q\ \left( |y| + R_\beta\ e^{-\tau} \right)^{q/(q-1)}\ \left\{ \left( 1 - e^{-q\tau} \right)^{-1/(q-1)} - 1 \right\} + \frac{q\ \gamma_q}{q-1}\ \left( |y| + |z|\ e^{-\tau} \right)^{1/(q-1)}\ |z|\ e^{-\tau} \\
& \le & \gamma_q\ \left( |y| + 1 \right)^{1/(q-1)}\ \left\{ \frac{q}{q-1} + |y| + 1 \right\}\ \left\{ \left( 1 - e^{-q\tau} \right)^{-1/(q-1)} - 1 + R_\beta\ e^{-\tau} \right\}
\eean
as $\tau\ge\log{R_\beta}$. Combining the above two estimates give
\bean
& & \left| h(\tau,y) - \sup_{z\in\RR^N}{\left\{ \left( h_0(z) - \gamma_q\ |y|^{q/(q-1)}\  \right)_+ \right\}} \right\vert\\
& \le & C(q)\ \left( |y| + 1 \right)^{q/(q-1)}\ \left\{ \left( 1 - e^{-q\tau} \right)^{-1/(q-1)} - 1 + R_\beta\ e^{-\tau} \right\} + 2\ \beta\,,
\eean
whence
\beqn
\label{ap5}
\left| h(\tau,y) - h_s(y) \right| \le C(q)\ \left( |y| + 1 \right)^{q/(q-1)}\ \left\{ \left( 1 - e^{-q\tau} \right)^{-1/(q-1)} - 1 + R_\beta\ e^{-\tau} \right\} + 2\ \beta
\eeqn
for $(\tau,y)\in [\log{R_\beta},\infty)\times\RR^N$. On the other hand, if $\tau\ge \log(R_\beta)$, $|y|\ge Y:=1+(\Vert h_0\Vert_\infty/\gamma_q)^{(q-1)/q}$ and $z\in\RR^N$, we have either $\left| y - z\ e^{-\tau} \right| \ge Y-1$ and 
\bean
& & h_0(z) - \gamma_q\ \left| y-z\ e^{-\tau} \right|^{q/(q-1)}\ \left( 1 - e^{-q\tau} \right)^{-1/(q-1)}\\ 
& \le & \left( 1 - e^{-q\tau} \right)^{-1/(q-1)}\ \left\{ \|h_0\|_\infty\ \left( 1 - e^{-q\tau} \right)^{1/(q-1)} - \gamma_q\ \left| y-z\ e^{-\tau} \right\vert^{q/(q-1)} \right\}\\
& \le & \left( 1 - e^{-q\tau} \right)^{-1/(q-1)}\ \left\{ \|h_0\|_\infty - \gamma_q\ (Y-1)^{q/(q-1)} \right\}\\
& \le & 0\,,
\eean
or $\left| y - z\ e^{-\tau} \right| < Y-1$ and 
$$
|z| \ge \left| y\ e^\tau \right| - \left| z - y\ e^\tau \right| \ge Y\ e^\tau - (Y-1)\ e^\tau = e^\tau \ge R_\beta\,,
$$
so that 
$$
h_0(z) - \gamma_q\ \left| y-z\ e^{-\tau} \right|^{q/(q-1)}\ \left( 1 - e^{-q\tau} \right)^{-1/(q-1)} \le \beta
$$
by \eqref{ap4}. Therefore,
\beqn
\label{ap6}
h(\tau,y) \le \beta \quad\mbox{ for }\quad (\tau,y)\in[\log{R_\beta},\infty)\times\RR^N\setminus B(0,Y)\,.
\eeqn

The claim \eqref{ap3} then easily follows from \eqref{ap5} and \eqref{ap6}. \qed



\begin{thebibliography}{99}  
 
\bibitem{BdCD97} 
Martino Bardi and Italo Capuzzo-Dolcetta, \textit{Optimal Control 
and Viscosity Solutions of Hamilton-Jacobi-Bellman Equations}, 
Systems Control Found. Appl., Birkh\"auser, Boston, 1997. 
 
\bibitem{Bl94} 
Guy Barles, \textit{Solutions de Viscosit\'e des Equations 
d'Hamilton-Jacobi}, Math\'ematiques \& Applications \textbf{17}, 
Springer-Verlag, Berlin, 1994. 
 
\bibitem{BP88} 
Guy Barles and Beno\^\i t Perthame, \textit{Exit time problems in 
optimal control and vanishing viscosity method}, SIAM J. Control 
Optim. \textbf{26} (1988), 1133--1148. 

\bibitem{BR06}
Guy Barles and Jean-Michel Roquejoffre, 
\textit{Ergodic type problems and large time behaviour of unbounded solutions of Hamilton-Jacobi equations},
Comm. Partial Differential Equations \textbf{31} (2006), 1209--1225. 

\bibitem{BS00}
Guy Barles and Panagiotis E. Souganidis,
\textit{On the large time behavior of solutions of Hamilton-Jacobi equations},
SIAM J. Math. Anal. \textbf{31} (2000), 925--939.
 
\bibitem{BtL08} Jean-Philippe Bartier and Philippe Lauren\c cot, 
\textit{Gradient estimates for a degenerate parabolic equation 
with gradient absorption and applications}, J. Funct. Anal. \textbf{254} (2008), 851--878. 

\bibitem{BKL04} 
Sa\"\i d Benachour, Grzegorz Karch, and Philippe Lauren\c cot, 
\textit{Asymptotic profiles of solutions to viscous 
Hamilton-Jacobi equations}, J. Math. Pures Appl. (9) \textbf{83} 
(2004), 1275--1308. 
 
\bibitem{CIL92} 
Michael G.~Crandall, Hitoshi Ishii, and Pierre-Louis Lions, 
\textit{User's guide to viscosity solutions of second order 
partial differential equations}, Bull. Amer. Math. Soc. (N.S.) 
\textbf{27} (1992), 1--67. 

\bibitem{EM94}
Juan R. Esteban and Pierangelo Marcati, 
\textit{Approximate solutions to first and second order quasilinear evolution equations via nonlinear viscosity}, Trans. Amer. Math. Soc. \textbf{342} (1994), 501--521.

\bibitem{Ev98} Lawrence C. Evans, Partial Differential Equations,
Grad. Stud. Math. \textbf{19}, Amer. Math. Soc., Providence, RI, 1998.

\bibitem{GGIS91} Yoshikazu Giga, Shun'ichi Goto, Hitoshi Ishii, and Moto-Hiko Sato, 
\textit{Comparison principle and convexity preserving properties for
singular degenerate parabolic equations on unbounded domains}, 
Indiana Univ. Math. J. \textbf{40} (1991), 443--470. 

\bibitem{Gi05} 
Brian H. Gilding, \textit{The Cauchy problem for $u_t=\Delta 
u+\vert \nabla u\vert^q$, large-time behaviour}, J. Math. Pures 
Appl. (9) \textbf{84} (2005), 753--785. 

\bibitem{GK04}
Brian H. Gilding and Robert Kersner, \textit{Travelling Waves in Nonlinear Diffusion-Convection Reaction}, Progr. Nonlinear Differential Equations Appl. \textbf{60}, Birkh\"auser, Basel, 2004. 

\bibitem{Ha93} Richard S. Hamilton, 
\textit{A matrix Harnack estimate for the heat equation}, Comm. Anal. Geom.  \textbf{1} (1993), 113--126.

\bibitem{I87} Hitoshi Ishii, 
\textit{A simple, direct proof of uniqueness for solutions of the 
Hamilton-Jacobi equations of eikonal type}, Proc. Amer. Math. Soc. 
\textbf{100} (1987), 247--251. 

\bibitem{I08} Hitoshi Ishii, 
\textit{Asymptotic solutions for large time of Hamilton-Jacobi equations in Euclidean $n$ space}, 
Ann. Inst. H. Poincar\'e Anal. Non Lin\'eaire \textbf{25} (2008), 231--266.

\bibitem{KV88} 
Shoshana Kamin and Juan Luis V\'azquez, 
\textit{Fundamental solutions and asymptotic behaviour for the $p$-Laplacian equation}, Rev. Mat. Iberoamericana \textbf{4} (1988), 339--354.

\bibitem{LS05} Philippe Lauren\c cot and Philippe Souplet, 
\textit{Optimal growth rates for a viscous Hamilton-Jacobi equation}, J. Evolution Equations \textbf{5} (2005), 123--135.

\bibitem{LV07} Philippe Lauren\c cot and Juan Luis V\'azquez, 
\textit{Localized non-diffusive asymptotic patterns for nonlinear parabolic equations with gradient absorption}, J. Dynamics Differential Equations \textbf{19} (2007), 985--1005.

\bibitem{LT01} Chi-Tien Lin and Eitan Tadmor, 
\textit{$L^1$-stability and error estimates for approximate Hamilton-Jacobi solutions}, Numer. Math. \textbf{87} (2001), 701--735.

\bibitem{NR99} 
Gawtum Namah and Jean-Michel Roquejoffre, \textit{Remarks on the 
long time behaviour of the solutions of Hamilton-Jacobi 
equations}, Comm. Partial Differential Equations \textbf{24} 
(1999), 883--893. 

\bibitem{Ro01} 
Jean-Michel Roquejoffre, \textit{Convergence to steady states or 
periodic solutions in a class of Hamilton-Jacobi equations}, J. 
Math. Pures Appl. (9) \textbf{80} (2001), 85--104. 
 
\bibitem{St02} 
Thomas Str\"omberg, 
\textit{The Hopf-Lax formula gives the unique viscosity solution},
Differential Integral Equations \textbf{15} (2002), 47--52.

\end{thebibliography}
\end{document}